\documentclass[a4paper]
{cedram-smai-jcm}

\usepackage[utf8]{inputenc}
\usepackage[british]{babel}
\usepackage{amsmath,amssymb}
\usepackage{tikz}
\usepackage{esint}

\newtheorem{Hypothesis}{Assumption}

\newcommand{\LLL}{\lvert\!\lvert\!\lvert}
\newcommand{\RRR}{\rvert\!\rvert\!\rvert}

\title[A posteriori estimates for biharmonic problems]
      {A posteriori error estimates for nonconforming discretizations
       of singularly perturbed biharmonic operators}

\author[D. Gallistl]{\firstname{Dietmar} \lastname{Gallistl}}
\address{Friedrich-Schiller-Universit\"at Jena, Institut f\"ur Mathematik,
         Ernst-Abbe-Platz 2, 07743 Jena, Germany}
\email{dietmar.gallistl@uni-jena.de}
\thanks{The first author is supported by the 
European Research Council
(ERC Starting Grant \emph{DAFNE}, agreement ID 891734).
The second author acknowledges support by 
Sino-German (CSC-DAAD) Postdoc Scholarship Program, 2021 ID 57575640%
}

\author[S. Tian]{\firstname{Shudan} \lastname{Tian}}
\address{School of Mathematics and Computational Science,
         Xiangtan University, Xiangtan 411105, China}
\email{shudan.tian@xtu.edu.cn}

\keywords{nonconforming finite element, singular perturbation,
          biharmonic, a~posteriori error estimation}
  
\subjclass{65N15; 65N30}

\begin{document}

\begin{abstract}
For the pure biharmonic equation and a biharmonic singular perturbation
problem, a residual-based error estimator is introduced which
applies to many existing nonconforming finite elements.
The error estimator involves the local best-approximation error
of the finite element function by piecewise polynomial functions
of the degree determining the expected approximation order,
which need not coincide with the maximal polynomial degree of
the element, for example if bubble functions are used.
The error estimator is shown to be reliable and locally efficient
up to this polynomial best-approximation error and oscillations
of the right-hand side.
\end{abstract}

\maketitle

\section{Introduction}

Given an open and bounded polytopal Lipschitz domain
$\Omega\subseteq \mathbb R^d$ in dimension $d\in\{2,3\}$, we
consider the following model problem of the biharmonic type.
Given a  parameter $\varepsilon\in (0,1]$ and $\alpha\in\{0,1\}$
and a right-hand side function $f$, we seek the solution $u$
to
\begin{equation}
	\begin{cases}
		\varepsilon^2\triangle^2 u-\alpha \triangle u = f\quad & \text{in }\Omega,  \\
		u = \frac{\partial u}{\partial n} =  0 \quad &\text{on }\partial\Omega.
	\end{cases}
	\label{e:modelproblem}
\end{equation}
Here, $n$ denotes the outer unit normal to $\partial\Omega$.
We focus on two choices of parameters.
For $\varepsilon=1$ and $\alpha=0$, we obtain the biharmonic equation
\begin{equation*}
	\triangle^2 u= f\quad  \text{in }\Omega
\end{equation*}
subject to clamped boundary conditions.
For the choice of a possibly small positive parameter $\varepsilon$
and $\alpha=1$, we obtain the singularly perturbed fourth-order problem
\begin{equation*}
	\varepsilon^2\triangle^2 u-\triangle u= f\quad  \text{in }\Omega
\end{equation*}
with clamped boundary conditions.
These are the two relevant scenarios because for $\alpha=0$, the problem
can always be reduced to the biharmonic equation after scaling of the 
right-hand side. If $\alpha=1$, the problem is singularly perturbed
and the solution $u_0$ to
\begin{equation*}
	-\triangle u_0 = f\quad  \text{in }\Omega
	\quad\text{and}\quad
	u_0  =  0 \quad \text{on }\partial\Omega
\end{equation*}
corresponds to the formal limit $\varepsilon=0$.
For convex domains, the convergence in the $H^1$ norm for $\varepsilon\to 0$
has been established in \cite[Lemma~5.1]{NTW}.
It is known that some nonconforming finite element methods are
convergent for the pure biharmonic equation but not for 
Poisson's equation and therefore perform poorly for the
singularly perturbed problem \eqref{e:modelproblem}.
For example, in \cite{NTW} the failure of the Morley
element was rigorously shown and a $C^0$ conforming alternative
similar to the Morley element was proposed.

For the biharmonic equation, a~posteriori error estimators were 
established for the Morley element \cite{BeiNiirSten2007}
and for some low-order rectangular elements \cite{CHDbipost}.
Under additional regularity assumptions, an a~posteriori error estimate
for the Specht element was derived in \cite{finitewang}. 
Beyond these references, there are no a~posteriori error estimates
available for nonconforming discretizations of \eqref{e:modelproblem}.
The difficulty typically encountered in the a~posteriori error analysis
of nonconforming finite element methods for the biharmonic equation
is that they usually do not contain a suitable conforming subspace.
In the reliability proof, the error can therefore not be approximated
by a conforming finite element function that is a valid test function
for the discrete problem. When the error is interpolated 
by a nonconforming function, care has to be taken that 
non-residual face terms resulting from an integration by parts
have suitable orthogonality properties.
Such properties fail to hold for most of the existing methods
because the maximal polynomial degree (henceforth denoted by $L$)
 of the shape functions is
higher than the degree of weak continuity across the elements
in the triangulation. Such high-degree polynomial are, however,
often relevant to the well-posedness (stability) of the method
(a set of shape functions for which the degrees of freedom are
linear independent is required), but 
irrelevant, in the generic case, to the approximation order
prescribed by the maximal degree $\ell$ of polynomials that are
fully included in the local shape function space.
The technique we propose in this work is therefore to include a
new term in the residual-based a~posteriori error estimator.
For a simplex $T$ in the triangulation and the finite element 
solution $u_h$ it reads
$$
\LLL u_h-\Pi_\ell u_h \RRR_{T}
$$
for the $L^2$ projection $\Pi_\ell$ to the polynomial
functions of degree not larger than $\ell$ over $T$.
The precise definition of the energy seminorm 
$\LLL \cdot \RRR_{T}$ can be found in
Section~\ref{ss:FEM}.
Scaling arguments and the triangle inequality
show that this term is locally
efficient up to the best approximation of the exact solution $u$
by piecewise polynomials of degree $\ell$.
It turns out that this approach enables a~posteriori error estimates
for a large class of nonconforming elements characterized by
the weak continuity from Assumption~\ref{assumption:main} below
and covers many existing methods for which no a~posteriori
estimates have been available before, see
Table~\ref{table:elementlist}.

An important application of this idea is the design of 
a~posteriori error estimators
for the fourth-order singular perturbation problem. 
For discontinuous $H^2$-nonconforming elements, the weak continuity 
properties of the function value jump are usually weaker than
those of the the gradient jump, see, e.g., \cite{Morley,Hu&Zhang2019}.
Consequently, such nonconforming elements tend to lose two orders of approximation
when the parameter $\varepsilon$ becomes small. 
Consequently, as mentioned above,
employing the Morley element directly in this problem does
not lead to convergence.
One remedy proposed in \cite{ModifyMorley} is to utilize a modified
bilinear form with the 
$P_1$ Lagrange interpolation in the second-order term.
This technique, known as the modified Morley method,
can also be extended to three-dimensional cases \cite{3DmodifiedMorley}.
The idea of the modified Morley element in 2D extends to other 
discontinuous $H^2$-nonconforming methods as outlined in \cite{Hu&Zhang2019}.
However such methods still do not reach the full convergence rate 
when $\varepsilon$ is small.    
Hence, $C^0$-conforming but $H^2$-nonconforming elements are particularly attractive
for the singular perturbation problem \cite{NTW,Guzm2012}.
In the literature, Wang and Zhang \cite{Wangsingular} developed 
error estimators for some low-order nonconforming elements assuming 
$H^3$ regularity.  Nevertheless, there have been no results under minimal 
regularity assumptions or for high-order methods or $d\geq3$ so far.
For the singularly perturbed Laplacian, parameter-robust a~posteriori error
estimators were derived \cite{Verfurth1998}. We combine these techniques
with the approach outlined above for proving a~posteriori error estimates
for the bi-Laplacian in the singularly perturbed case. 
As an alternative approach to nonconforming schemes (not included in our 
analysis),
$C^0$-IPDG methods were proposed for fourth order problems \cite{C0IPDG},
and a~posteriori techniques were considered in \cite{BrennerGudiSung2010}
and in \cite{PostC0IPDGS} for the singularly perturbed situation.

The remaining parts of this article are organized as follows.
Section~\ref{s:prelim} provides preliminaries on the model problem
and the finite element discretization. The reliability and efficiency
of the error estimator are proven in Section~\ref{s:estimator}.
Section~\ref{s:examples} is devoted to examples and Section~\ref{s:num}
reports numerical experiments. For convenient reading, some technical details
are listed as Appendix~\ref{ax:HCT} and \ref{ax:localizationproof}.

\bigskip
Standard notation for function spaces is used throughout
the article. The inner product in the Lebesgue space
$L^2(\omega)$ with respect to $\omega\subseteq \Omega$
is denoted by $(\cdot,\cdot)_{0,\omega}$, and the $L^2$
norm over $\omega$ is denoted by $\|\cdot\|_{0,\omega}$.
If $\omega=\Omega$, we simply write 
$(\cdot,\cdot)_0$ and $\|\cdot\|_0$.
The second-order $L^2$-based Sobolev space with clamped boundary
condition is denoted by $V:=H^2_0(\Omega)$.
The Euclidean inner product of $x,y\in\mathbb R^n$
is denoted by $x\cdot y$.
The divergence operator for matrix-valued functions is understood
row-wise.
The notation $a\lesssim b$ means the inequality
$a\leq C b$ up to a multiplicative constant $C$ that does not depend
on the mesh-size $h$ or the parameter $\varepsilon$.

\section{Preliminaries}\label{s:prelim}

\subsection{Weak formulation of the model problem}

The weak formulation of \eqref{e:modelproblem} is based on the Sobolev space
$V:=H^2_0(\Omega)$.
Given $f\in L^2(\Omega)$, it seeks $u\in V$ such that
\begin{equation}\label{e:weak_modelproblem}
a(u,v) 
= (f,v)_0
\quad \text{for all } v\in V
\end{equation}
for the bilinear form
$$
a(u,v) = \varepsilon^2(D^2u,D^2v)_0 + \alpha (\nabla u,\nabla v)_0.
$$
We denote the induced energy norm by $\|v\|_a:=a(v,v)^{1/2}$.

\subsection{Discrete functions}
Throughout this work, we assume an underlying regular 
simplicial triangulation
$\mathcal T$ of the domain $\Omega$ from a shape-regular family.
For each element $T\in\mathcal{T}$, the diameter of $T$ is denoted
by $h_T$, and by $h$ we denote the piecewise constant mesh-size function
defined by $h|_T=h_T$ for any $T\in \mathcal T$.
The set of $(d-1)$-dimensional hyperfaces is denoted by $\mathcal{F}$.
The subset of interior faces is denoted by $\mathcal{F}(\Omega)$.
The diameter of a face $F$ is denoted by $h_F$.
Any interior face $F\in\mathcal F(\Omega)$ is shared by two elements, which we
denote by $T_+,T_-\in\mathcal T$. If $F$ is a boundary face, $T_+$ is
the unique element containing $F$ and we use the convention $T_-=\emptyset$.
The interior of $T_+\cup T_-$ is denoted by $\omega_F$ and referred to
as the face patch of $F$.
To every face $F\in\mathcal F$ we assign a unit normal vector $n_F$.
If $F\subseteq \partial \Omega$ is a boundary face, we choose $n_F$ to be
outward-pointing with respect to $\Omega$.
The set of vertices in the triangulation is denoted by $\mathcal N$.
Given a vertex $y\in\mathcal N$, its vertex patch $\omega_y$ is the 
interior of the union of all elements of $\mathcal T$ that contain
$y$; the diameter of the vertex patch is denoted by $h_y$.
The jump and the average of a function $v$ on an 
interior face $F\in \mathcal F(\Omega)$ are defined as
$[v]|_F: = v_+ - v_-$ and $\{v\}|_F:=2^{-1}(v_+ + v_-)$.
Here we set $v_\pm:= v|_{T_\pm}$.
If $F\subseteq\partial\Omega$ is a boundary face, we 
set $[v]|_F=\{v\}_F=v|_{T_+}$.

The symbols $\nabla_h$, $\operatorname{div}_h$, $D^2_h$, $\triangle_h$ denote the
$\mathcal T$-piecewise
application of the differential operators 
$\nabla$, $\operatorname{div}$, $D^2$, $\triangle$,
respectively.
Given a subdomain $\omega\subseteq\bar\Omega$, the space of polynomial
functions over $\omega$ of degree not larger than $\ell\geq 0$ is
denoted by $P_\ell(\omega)$. The piecewise polynomial functions
with respect to $\mathcal T$ are denoted by
$$
P_\ell(\mathcal T) = \{v\in L^2(\Omega): v|_T\in P_\ell(T) 
                         \text{ for any }T\in\mathcal T\}.
$$
The $L^2(\Omega)$-orthogonal projection onto $P_\ell(\mathcal T)$
is denoted by $\Pi_\ell$.

\subsection{Finite element discretization}\label{ss:FEM}

The finite element space
$V_h$ is assumed to be a subspace of the piecewise polynomial functions
with some given maximal degree $L$, i.e., $V_h\subseteq P_{L}(\mathcal T)$.
On $V+V_h$, we introduce the semidefinite bilinear form
$$
a_h(u,v) = \varepsilon^2(D_h^2u,D_h^2v)_0 + \alpha (\nabla_h u,\nabla_h v)_0
$$
with seminorm $\|v\|_{a,h}:=a_h(v,v)^{1/2}$.
We assume throughout that any element $v_h\in V_h$ satisfies the weak
continuity condition
\begin{equation}
	\label{e:weak_cont}
	\|[v_h]_F\|_{0,F}
	\lesssim
	h_F
	\|\nabla_F [v_h]_F\|_{0,F}
	\quad\text{for all }F\in\mathcal F
\end{equation}
where $\nabla_F$ is the tangential gradient operator on $F$.

\begin{remark}
 We refer to condition \eqref{e:weak_cont} as weak continuity because it obviously
 excludes constant jumps across the faces.
\end{remark}

Condition \eqref{e:weak_cont} and Assumption~\ref{assumption:main}
below guarantee that $a_h$ is positive definite over $V_h$
such that $\|\cdot\|_{a,h}$ is a norm on $V_h$,
see Lemma~\ref{l:posdef} for a proof.
We further assume the existence of an operator
$J_h:V\to V_h$ with the approximation and stability properties
\begin{subequations}\label{e:quasiint}
\begin{equation}\label{e:quasiint_a}
  \|h^{-2} (v-J_h v)\|_0
  +\|h^{-1} \nabla_h (v-J_h v)\|_0
  +\|D_h^2 (v-J_h v)\|_0
  \lesssim
  \|D^2 v\|_0
\end{equation}
and
\begin{equation}\label{e:quasiint_b}
  \|h^{-1} (v-J_h v)\|_0
  +\|\nabla_h (v-J_hv)\|_0
  \lesssim
  \|\nabla v\|_0
\end{equation}
\end{subequations}
for all $v\in V$.
Such operators, referred to as quasi-interpolation operators,
are usually defined as a regularization by local averaging,
see \cite{ErnGuermond_FEM1}.
We will also use,
for measurable subsets $\omega\subseteq\Omega$,
 the localized seminorm
$$
\LLL w \RRR_{\omega}
:=
\left(\varepsilon^2\|D^2_h w\|_{0,\omega}^2 + \alpha\|\nabla_h w\|_{0,\omega}^2\right)^{1/2}
\quad\text{for any } w\in V+V_h,
$$
for which we have 
$\LLL w \RRR_{\Omega} = \|w\|_{a,h}$.
An application of inverse inequalities shows that the $L^2$ projection
$\Pi_\ell$ is quasi-optimal with respect to the seminorm
$\LLL \cdot \RRR_{T}$, i.e.,
\begin{equation}\label{e:L2quasiopt}
 \LLL w_L-\Pi_\ell w_L\RRR_{T}
 \lesssim
 \min_{v_\ell\in P_\ell(T)}
 \LLL w_L-v_\ell\RRR_{T}
 \quad\text{for any }w_L\in P_L(T).
\end{equation}
While $L$ denotes the highest polynomial degree used for the 
local shape functions of $V_h$, usually $P_L(T)$ will not be
a subset of the local shape function space defining the underlying
finite element; a typical instance of such incomplete polynomial
shape function spaces is the enrichment by bubble functions.
It is instead assumed that the shape function space contains
$P_\ell(T)$ for some $\ell\leq L$.
Following \cite{Hu&Zhang2019},
we formulate the following consistency condition for the inter-element
jumps.

\begin{Hypothesis}\label{assumption:main}
There exists $\ell\geq2$ such that
any element $v_h\in V_h$ from the finite element space $V_h$ with respect to
$\mathcal{T}$ the following is satisfied.
\begin{itemize}
\item
If $\alpha=0$,
for any $F\in\mathcal F$, there holds
\begin{subequations}\label{eq:assumption_main}
	\begin{align}
	\label{eq:assumption_main_part1}
	& \int_F [ v_h]_F q\mathrm{d}S= 0 
	&
	&\text{for all } q\in 
        P_{\ell-3}(F),
	\\
	\label{eq:assumption_main_part2}
	& \int_F [\nabla_h v_h]_F\cdot Q\mathrm{d}S =0
	&
	&\text{for all } Q\in (P_{\ell-2}(F))^d .
	\end{align}
\end{subequations}
Here, we follow the convention that $P_{-1}(F)=\{0\}$.
\item
If $\alpha=1$,
we assume the inclusion $V_h\subseteq H_0^1(\Omega)$
and that \eqref{eq:assumption_main_part2} holds for all 
$F\in\mathcal F$.
\end{itemize}
\end{Hypothesis}

\begin{remark}
 Recall that in this work we only consider $\alpha\in\{0,1\}$.
 Clearly, the requirements for the case $\alpha=1$ in 
 Assumption~\ref{assumption:main} imply those for 
 $\alpha=0$. The additional assumption for $\alpha=1$ is that the 
 weak continuity \eqref{eq:assumption_main_part1} is replaced
 by strong continuity.
\end{remark}

The finite element discretization of \eqref{e:weak_modelproblem} seeks
$u_h\in V_h$ such that
\begin{equation}\label{e:discr_modelproblem}
	a_h(u_h,v_h)=(f,v_h)_0
	\quad\text{for all }v_h\in V_h.
\end{equation}
The well-posedness of the discrete problem is a consequence of the following
lemma.
\begin{lemma}\label{l:posdef}
 The weak continuity \eqref{e:weak_cont} and Assumption~\ref{assumption:main}
 imply that $a_h$ is positive definite over $V+V_h$.
\end{lemma}
\begin{proof}
 Any $v\in V+V_h$ with $a_h(v,v)=0$ is piecewise affine.
 The combination of \eqref{e:weak_cont} and \eqref{eq:assumption_main_part2}
 (recall $\ell\geq 2$) therefore shows that $v$ is globally continuous.
 Since by \eqref{eq:assumption_main_part2} the gradient does not jump across
 the inter-element faces, $v$ is globally affine.
 By employing \eqref{e:weak_cont} and \eqref{eq:assumption_main_part2} on the 
 boundary faces, we deduce $v=0$.
\end{proof}

\section{A posteriori error estimates}\label{s:estimator}

\subsection{Preparatory identities}

We use the following error decomposition \cite{CHDbipost}
based on the Pythagorean identity.  
\begin{lemma}\label{l:errordecomposition}
Any $u\in V$ and $u_h\in V_h$ satisfy
\begin{equation*}
	\|u-u_h\|_{a,h}^2
	= \sup_{v\in V,\|v\|_a=1 }  |a_h(u-u_h,v)|^2
	 +\min_{v\in V}\|u_h-v\|_{a,h}^2.	
\end{equation*}
\end{lemma}

In order to estimate the first term on the right-hand side of 
the foregoing error decomposition, we note some basic identities
based on integration by parts.

\begin{lemma}[integration by parts]\label{l:int-by-parts}
 Let Assumption~\ref{assumption:main} hold.
 Let $u_h\in V_h$ and $w\in V+V_h$ be given.
 Then,
 for $\ell$ as in Assumption~\ref{assumption:main}
 and any $p_h\in P_\ell(\mathcal T)$,
 $$
   (D_h^2 u_h, D_h^2 w)_0
   =
   (\triangle_h^2 u_h, w)_0
   +
   T_1
   +T_2
$$
where
$$
 T_1:=
   \sum_{F\in \mathcal{F}(\Omega)}
      \left(\int_F[D_h^2u_h ]_Fn_F\cdot\{\nabla_h w\}_F\mathrm{d}S
   -
     \int_F[\operatorname{div}_hD_h^2u_h]_F \cdot n_F\{w\}_F\mathrm{d}S
     \right)
$$
and
$$
 T_2:=
   \sum_{F\in\mathcal F}
   \left(
   \int_F \{D_h^2(u_h-p_h)\}_F n_F\cdot[\nabla_h w]_F\mathrm{d}S 
   - \int_F \{\operatorname{div}_hD_h^2(u_h-p_h)\}_F\cdot n_F[w]_F\mathrm{d}S 
   \right).
 $$
 If $\alpha=1$, we have
 $$
   (\nabla_h u_h, \nabla_h w)_0
   =
   -(\triangle_h u_h, w)_0
   +
   T_3
$$
where
$$
 T_3:=
   \sum_{F\in \mathcal{F}(\Omega)}\int_F[\nabla_h u_h]_F\cdot n_F\{w\}_F\mathrm{d}S.
 $$
\end{lemma}
\begin{proof}
 The first identity follows from piecewise integration by parts and the 
 properties \eqref{eq:assumption_main}, which guarantee that 
 on every face $F$, the functions
 $[\nabla_h w]_F$ and $[w]_F$ are orthogonal to $\{D_h^2 p_h\}_F n_F$
 and $\{\operatorname{div}_h D_h^2 p_h\}_F\cdot n_F$, respectively. 
 The second stated identity follows from piecewise integration by parts and the 
 inclusion $V_h\subseteq H^1_0(\Omega)$
 guaranteed by Assumption~\ref{assumption:main} for the case $\alpha=1$.
\end{proof}

\subsection{Error estimator contributions}

As in \cite{Verfurth1998}, we define
the comparison coefficient $\kappa_T=\min\{1,h_T/\varepsilon\}$
of the mesh size $h_T$ and $\varepsilon$.
Similarly, we write
$\kappa_F=\min\{1,h_F/\varepsilon\}$
and
$\kappa_y=\min\{1,h_y/\varepsilon\}$
for a face $F$ or a vertex $y$.
We denote
\begin{align*}
    &\mu_0^2(T):= \LLL  u_h-\Pi_\ell u_h\RRR_{T}^2,
    &\;
	&\mu_1^2(T)
	  :=\kappa_T^2h_T^2
	   \|\varepsilon^2\triangle_h^2 u_h-\alpha \triangle_hu_h-f\|^2_{0,T},
    \\
	&\mu_2^2(F):= \varepsilon^3\kappa_F\|[D_h^2u_h]_F n_F\|_{0,F}^2,	
	&
	&\mu_3^2(F):= \kappa^2_F h_F 
	 \| [\alpha\nabla_h u_h-\varepsilon^2\mathrm{div}_h D^2_hu_h]_F\cdot n_F\|^2_{0,F}
	&
	&
\end{align*}
for any $T\in\mathcal T$ and any interior face $F\in\mathcal F(\Omega)$.
For boundary faces $F\subseteq\partial\Omega$
 we set $\mu_2^2(F)=\mu_3^2(F)=0$.
We abbreviate
$$
\mu_j(\mathcal T):=\sqrt{\sum_{T\in \mathcal T}\mu_j^2(T)},
\quad
\mu_k(\mathcal F):=\sqrt{\sum_{F\in \mathcal F}\mu_k^2(F)},
\quad
\text{ for } j=0,1 \text{ and }k=2,3.
$$
We further define, for any face $F\in\mathcal F$,
\begin{align*}
 \xi^2(F):=
\frac{\varepsilon}{\kappa_F} \|[\nabla_h u_h]_F\cdot n_F\|^2_{0,F}
	+  \frac{1}{\varepsilon\kappa_F^3} \|[u_h]_F\|^2_{0,F} 
    \;\text{and}\;
    \xi(\mathcal F) := \sqrt{\sum_{F\in\mathcal F} \xi^2(F)}.
\end{align*}

\subsection{Error estimator reliability}

We begin by bounding the first term of the decomposition from
Lemma~\ref{l:errordecomposition}.

\begin{lemma}[equilibrium residual]\label{l:equilibrium}
	Let Assumption~\ref{assumption:main} hold.
	Let $u\in V$ solve \eqref{e:weak_modelproblem} and let $u_h\in V_h$
	solve \eqref{e:discr_modelproblem}. Then, for any $v\in V$
	with $\|v\|_a=1$ we have
	\begin{align*}
    a_h(u-u_h,v)
	& 
	\lesssim 
	\mu_0(\mathcal T) + \mu_1(\mathcal T)
	+
	\mu_2(\mathcal F) + \mu_3(\mathcal F) .
\end{align*}
\end{lemma}
\begin{proof}
We consider the quasi-interpolation $J_h v\in V_h$ of the function $v$
satisfying \eqref{e:quasiint} and 
abbreviate $w:= v-J_h v$.
We deduce from the solution properties \eqref{e:weak_modelproblem}
and \eqref{e:discr_modelproblem}
of $u$ and $u_h$ that
\begin{equation*} 
	a_h(u-u_h,v) = (f,w)_0 - a_h(u_h, w).
\end{equation*}
We apply integration by parts to the second term on the right-hand side:
from Lemma~\ref{l:int-by-parts} with
$p_h := \Pi_\ell u_h$ we obtain
$$
(f,w)_0- a_h(u_h, w)
= (f-\varepsilon^2\triangle_h^2 u_h+\alpha\triangle_h u_h,w)_0
- \varepsilon^2 T_1 - \varepsilon^2 T_2 
- \alpha T_3
$$
with $T_1, T_2, T_3$ as in 
Lemma~\ref{l:int-by-parts}.
For the first term on the right-hand side of this identity,
\eqref{e:quasiint} leads to
$$
(f-\varepsilon^2\triangle_h^2 u_h+\alpha\triangle_h u_h,w)_0
\lesssim
\mu_1(\mathcal T).
$$
The trace inequality, properties \eqref{e:quasiint},
and the bounded overlap of patches show
$$
\sum_{F\in\mathcal F}
\kappa_F^{-2}h_F^{-1}\|\{w\}_F\|_{0,F}^2
\lesssim
\sum_{F\in\mathcal F}
\kappa_F^{-2} (h_F^{-2} \|w\|_{0,\omega_F}^2 + \|\nabla_h w\|_{0,\omega_F}^2)
\lesssim
\| v \|_a=1.
$$
Similarly, with the multiplicative trace inequality
$$
\|\{\nabla_h w\}_F\|_{0,F}
\lesssim
h_F^{-1/2} \|\nabla_h w\|_{0,\omega_F}
+
\|\nabla_h w\|_{0,\omega_F}^{1/2} \|D^2_h w\|_{0,\omega_F}^{1/2}
$$
we obtain
$$
\sum_{F\in\mathcal F}\varepsilon\kappa_F^{-1} \|\{\nabla w\}_F\|_{0,F}^2
\lesssim
\| v \|_a=1.
$$
The Cauchy inequality and the foregoing trace estimates therefore
establish
$$
|\varepsilon^2 T_1+\alpha T_3|
\lesssim \mu_2(\mathcal F) + \mu_3(\mathcal F).
$$
Similarly, trace and inverse inequalities and \eqref{e:quasiint}
give
$$
|\varepsilon^2 T_2| 
\lesssim
\mu_0(\mathcal T).
$$
This concludes the proof.
\end{proof}

We proceed with estimating the second term in the error decomposition
of Lemma~\ref{l:errordecomposition}, which refer to as the conformity
residual. 
As in \cite{CHDbipost}, we shall bound that term from above by designing 
a suitable conforming finite element approximation to $u_h$.
The construction is based on averaging operators employing generalized
Hsieh--Clough--Tocher (HCT) spaces, 
the details of which are summarized in Appendix~\ref{ax:HCT}
for convenient reading.
The general design of such averaged approximations is well known,
see \cite{BrennerGudiSung2010,HighorderHCT}.
In the case of a singular perturbation, such average approximation may
be required on some local sub-mesh.
We use the following localization argument.

\begin{lemma}[localization]\label{l:localization}
 Given any $u_h\in V_h$, we have
 \begin{align*}
    &\min_{v\in V} \| u_h - v \|_{a,h}^2
    \\
    &\lesssim
    \sum_{y\in\mathcal N}
    \min_{v_y\in V(\omega_y)} 
    \left(
    \frac{1}{h_y^2\kappa_y^2} \| u_h - v_y\|_{0,\omega_y}^2
    + \frac{1}{\kappa_y^2} \| \nabla_h(u_h - v_y)\|_{0,\omega_y}^2
    + \varepsilon^2 \| D^2_h(u_h - v_y)\|_{0,\omega_y}^2
    \right)
 \end{align*}
where $V(\omega_y)$ is the space of functions over 
the vertex patch $\omega_y$ 
that admit an extension to $\Omega$ that belongs to $V$.
\end{lemma}
\begin{proof}
 The proof is postponed to Appendix~\ref{ax:localizationproof}.
\end{proof}

\begin{lemma}[conformity residual]\label{l:conformity}
Let Assumption~\ref{assumption:main} hold.
	For $u_h\in V_h$, we have
	    $
		\min_{v\in V}\|u_h-v\|_{a,h}
		\lesssim 
		\xi(\mathcal F).
		$
\end{lemma}
\begin{proof}
 We begin with the localization from Lemma~\ref{l:localization}
 and consider a vertex $y\in\mathcal N$ with vertex
 patch $\omega_y$ of diameter $h_y$.
 We denote by $\mathcal F(y)$ the set faces $F\in\mathcal F$
 with $y\in F$.
 We denote by $\hat {\mathcal  T}=\hat {\mathcal  T}(\omega_y)$
 a uniformly refined triangulation
 of the patch $\omega_y$ of mesh size 
 \begin{equation}\label{e:hatscaling}
  \max_{\hat T\in\hat{\mathcal{T}}} \operatorname{diam} (\hat T)
  =:
  \hat h \leq \varepsilon \lesssim \kappa_y^{-1} \hat h
 \end{equation}
 Let $\hat V_C\subseteq V(\omega_y)$ denote
 the conforming generalized
 HCT finite element space 
 \cite{Ciarlet2002,HighorderHCT,WSNeilan}
 whose local shape functions
 contain $P_L(\hat T)$ for any $\hat T\in\hat{\mathcal T}$.
 The operator $\hat \Pi_C:V_h|_{\omega_F}\to \hat V_C$
 assigns a lifted object $\hat \Pi_C u_h $
 to $u_h|_{\omega_z}$ by setting any global degree of freedom 
 (with respect to $\hat V_C$)
 as the average of the evaluation of the local node functionals applied to
 $u_h|_{\omega_z}$.
 Since the degrees of freedom of the HCT functions are point evaluations
 of the function and the gradient in the vertices, evaluations of the normal
 derivative over the faces, and evaluations inside the element,
 standard arguments with triangle and inverse inequalities show that
 \begin{align*}
  \hat h^{-4} \| u_h - \hat \Pi_C u_h\|_{0,\hat T}^2
  \lesssim
  \sum_{\substack{\hat F\in \hat {\mathcal{F}}(\hat T) \\ 
                  \exists F\in\mathcal F: 
                  \hat F\subseteq F \in\mathcal F(y)
                  }
        }
     (\hat h^{-1} 
	    \|[\nabla_h u_h]_{\hat F}\cdot n_{\hat F}\|^2_{0,\hat F}
	    +  \hat h^{-3} \|[u_h]_{\hat F}\|^2_{0,\hat F} )
 \end{align*}
 for any $\hat T\in\hat{\mathcal T}$ with set of faces 
 $\hat{\mathcal F}(\hat T)$.
 The sum is taken over all (fine) faces of $\hat T$ that are included
 in an element of $\mathcal F(y)$.
 Since $u_h$ is a polynomial function in each coarse simplex,
 its jumps are only nontrivial along the faces $F\in\mathcal F(y)$. We conclude
  \begin{align*}
  \hat h^{-4} \| u_h - \hat \Pi_C u_h\|_{0,\omega_y}^2
  \lesssim
  \sum_{F\in\mathcal F(y)}
              (\hat h^{-1} \|[\nabla_h u_h]_F\cdot n_F\|^2_{0,F}
                         +  \hat h^{-3} \|[u_h]_F\|^2_{0,F} )
  .
 \end{align*}
 Using inverse estimates over $\hat{\mathcal{T}}$
 and the above scaling \eqref{e:hatscaling} of $\hat h$,
 we thus obtain
 \begin{align*}
 &
    \frac{1}{h_y^2\kappa_y^2} \| u_h - \hat \Pi_C u_h\|_{0,\omega_y}^2
    + \frac{1}{\kappa_y^2} \| \nabla_h(u_h - \hat \Pi_C u_h)\|_{0,\omega_y}^2
    + \varepsilon^2 \| D^2_h(u_h - \hat \Pi_C u_h)\|_{0,\omega_y}^2
 \\
 &
 \qquad
  \lesssim
  \sum_{F\in\mathcal F(y)}
    \left[
    \frac{\varepsilon}{\kappa_F} \|[\nabla_h u_h]_F\cdot n_F\|^2_{0,F}
	+  \frac{1}{\varepsilon\kappa_F^3} \|[u_h]_F\|^2_{0,F} 
	\right]
  = \sum_{F\in\mathcal F(y)}\xi^2(F)
  .
 \end{align*}
Inserting this upper bound in the sum on the right hand side
of the displayed formula from Lemma~\ref{l:localization}
concludes the proof.
\end{proof}

We conclude the reliability of the error estimator
$$
\eta(\mathcal T) :
 = \mu_0(\mathcal T) +  \mu_1(\mathcal T)
  +  \mu_2(\mathcal F) +  \mu_3(\mathcal F) 
  + \xi(\mathcal F). 
$$

\begin{theorem}[reliability]
 Let Assumption~\ref{assumption:main} hold.
 The solutions $u\in V$ to \eqref{e:weak_modelproblem} 
 and $u_h\in V_h$ to \eqref{e:discr_modelproblem} satisfy
 $$
   \| u-u_h\|_{a,h} \lesssim \eta(\mathcal{T}).
 $$
\end{theorem}
\begin{proof}
 The result follows from combining
 Lemma~\ref{l:equilibrium} and Lemma~\ref{l:conformity}.
\end{proof}

\subsection{Error estimator efficiency}
The error estimator $\eta$ is locally efficient up to data
oscillations and the best-approximation by piecewise
polynomials of degree $\ell$.
The proofs partly
follow standard arguments \cite{VRposteriori,Verfurth1998} and we only
highlight some important aspects in the analysis.
Given $f\in L^2(\Omega)$, we define its oscillations of
order $\ell$ with respect to $T\in\mathcal T$ as
$$
\operatorname{osc}^2(f,\ell,T) = 
h_T^2\kappa_T^2 \|f-\Pi_\ell f\|^2_{0,T}.
$$

\begin{theorem}[local efficiency]
 Let Assumption~\ref{assumption:main} hold.
 Let $u\in V$ solve \eqref{e:weak_modelproblem} 
 and $u_h\in V_h$ solve \eqref{e:discr_modelproblem}.
 The local error estimator contributions satisfy for any $T\in\mathcal T$
 and any $F\in\mathcal F$ that
 \begin{align*}
  \xi^2(F)
   &\lesssim \LLL  u-u_h \RRR_{\omega_F}^2,
 \\      
 \mu_0^2(T)
 &\lesssim \LLL  u-u_h \RRR_{T}^2
 +
 \min_{p_h\in P_\ell(T)}
         \LLL  u-p_h\RRR_{T}^2,
 \\
 \mu_1^2(T) 
   &\lesssim \LLL  u-u_h \RRR_{T}^2 + \operatorname{osc}^2(f,\ell,T),
 \\
 \mu_2^2(F) + \mu_3^2(F) 
   &\lesssim \LLL  u-u_h \RRR_{\omega_F}^2
      + \sum_{K\in \{T_+,T_-\}}\operatorname{osc}^2(f,\ell,K)
      .
 \end{align*}
\end{theorem}

\begin{proof}
We begin by estimating $\xi^2(F)$.
If $\alpha=0$, the efficiency follows from standard arguments
as in \cite{BeiNiirSten2007,HuShi2009}.
The same applies to the case $\alpha=1$ and $h_F\leq\varepsilon$.
Let now $\alpha=1$ and $\varepsilon\leq h_F$.
In particular, $u_h$ is continuous due to Assumption~\ref{assumption:main}.
We have $\kappa_F=1$
and therefore
the weighted trace inequality 
\begin{align*}
\varepsilon \|[\nabla_h u_h]_F\cdot n_F\|^2_{0,F}
\lesssim
\frac{\varepsilon}{h_F}
\|\nabla_h (u-u_h)\|_{0,\omega_F}^2
+
\varepsilon
\|\nabla_h (u-u_h)\|_{0,\omega_F}
\|D^2_h (u-u_h)\|_{0,\omega_F}
\end{align*}
proves the bound
\begin{align*}
\xi^2(F)
\lesssim \LLL u-u_h\RRR_{\omega_F}^2.
\end{align*}

The local efficiency of $\mu_0^2(T)$ follows from the triangle inequality
and the quasi-optimality of the piecewise $L^2$ projection
displayed in \eqref{e:L2quasiopt}.
The local efficiency of $\mu_1^2(T)$ follows from a rather
standard argument \cite{VRposteriori,Verfurth1998}.
We include a brief proof to highlight the role
of the parameter $\kappa_T$.
We denote by $b_T=\prod_{j=1}^{d+1}\lambda_j$ (with barycentric coordinates
$\lambda_j$ over $T$) the usual volume bubble function.
We define $\psi_T:=b_T^2(\varepsilon^2\triangle_h^2 u_h-\triangle_h u_h-\bar f)$
for $\bar f:=\Pi_\ell f$.
Equivalence of finite-dimensional norms and the scaling 
$\|b_T^2\|_{L^\infty(T)}\leq 1$
show
$$
\mu_1^2(T) \lesssim \kappa_T^2 h_T^2 
    \left[ (\varepsilon^2\triangle_h^2 u_h-\alpha \triangle_hu_h- f,\psi_T)_0
     + \| f-\bar f\|_{0,T}^2 \right]
$$
so that, after integration by parts, we conclude
$$
\mu_1^2(T) \lesssim \kappa_T^2 h_T^2 
      \LLL u-u_h\RRR_{T}
           \LLL \psi\RRR_{T}
     + \kappa_T^2 h_T^2 \| f-\bar f\|_{0,T}^2  .
$$
The inverse estimate and the pointwise bound on $b_T$ imply
$$
\LLL \psi\RRR_{T}
\lesssim
h_T^{-1}\kappa_T^{-1} \|\varepsilon^2\triangle_h^2 u_h-\alpha \triangle_hu_h- \bar f\|_{0,T}
$$
such that
$$
\mu_1^2(T) \lesssim 
   \LLL u-u_h\RRR_{T}^2
     + \kappa_T^2 h_T^2 \| f-\bar f\|_{0,T}^2 .
$$

We proceed with the proof of efficiency of $\mu_2^2(F)$
for an interior face $F\in\mathcal F(\Omega)$.
Since the norm tangential-normal jump of $D^2_h u_h$ can be bounded
by the terms included in $\xi^2(F)$ after an inverse estimate along
$F$, we focus on bounding the normal-normal jump of the discrete Hessian.
To this end, we follow an idea from
\cite{Verfurth1998} and use a face bubble function 
$$
 \chi_{F,\delta}:= (b_{T_+,\delta}- \frac{|T_-|}{|T_+|}b_{T_-,\delta})b_{F,\delta}
$$
where the volume and face bubble functions
$$
b_{K,\delta}:=\prod_{j=1}^{d+1}\lambda_{j,K,\delta}
\quad\text{and}\quad
b_{F,\delta}:=\prod_{\substack{j=1 \\ z_j\in F}}^{d+1}\lambda_{j,K,\delta}
$$
for a simplex $K$ with vertices $z_1,\dots,z_{d+1}$ and a face $F\subseteq K$
are defined with respect to the barycentric coordinates $\lambda_{j,K,\delta}$
of a subsimplex of $K$ that contains $F$ and has a height over
$F$ of order $\delta h_K$.
The product in the definition of $b_{F,\delta}$ runs over the $d$ vertices
satisfying $z_j\in F$.
The bubble functions are extended by zero
to the patch $\overline{\omega}_F$. Details on the construction can be found in
\cite{Verfurth1998}.
The function $\chi_{F,\delta}$ vanishes on $F$ and belongs to
$H^2_0(\omega_F)$.
The bounds 
$$
  \|\chi_{F,\delta}\|_{0,\omega_F}\lesssim (h_F\delta)^{1/2} |F|^{1/2}
 \quad\text{and}\quad
 \|\frac{\partial \chi_{F,\delta}}{\partial n_F}\|_{L^\infty(F)}\lesssim (h_F\delta)^{-1}
$$ 
can be verified by scaling arguments with the shape regularity.
The function
$$
\psi:= \chi_{F,\delta} n_F^\top [D_h^2 u_h]_F n_F
$$
therefore satisfies, after a continuation to the patch $\omega_F$
as in \cite{VRposteriori}, that
$$
 |\psi|_{m,\omega} \lesssim (\delta h)^{1/2-m}
         \| n_F^\top [D_h^2 u_h]_F n_F\|_{0,F},
         \quad m=0,1,2.
$$

With equivalence of norms and integration by parts we deduce
\begin{align*}
\varepsilon^3 \kappa_F \| n_F^\top [D_h^2 u_h]_F n_F\|_{0,F}^2
&\lesssim
\varepsilon^3 \kappa_F h_F \delta 
   ([D_h^2 u_h]_F n_F,\nabla\psi)_{0,F}
\\
&=
\varepsilon \kappa_F h_F \delta 
(\varepsilon^2 (D^2_h u_h,D^2\psi)_0 - (\varepsilon^2 \Delta_h^2 u_h,\psi)_0 ).
\end{align*}
After adding and subtracting $\Delta_h u_h$ and $f$ we therefore obtain
\begin{align*}
&\varepsilon^3 \kappa_F \| n_F^\top [D_h^2 u_h]_F n_F\|_{0,F}^2
\\
&\lesssim
\varepsilon \kappa_F h_F \delta 
(a_h( u_h,\psi) - (f,\psi)_0) - (\varepsilon^2 \Delta_h^2 u_h-\Delta_h u_h-f,\psi)_0 )
\\
&
\lesssim
\varepsilon \kappa_F h_F \delta
(\LLL u-u_h\RRR_{\omega_F} \|\psi\|_a
 + (\mu_1(T_+)+\mu_1(T_-))  \kappa_F^{-1}h_F^{-1}\|\psi\|_0)
.
\end{align*}
The choice $\delta:=\min\{1,\varepsilon/h_F\}$ and direct
computations lead to
$$
\varepsilon \kappa_F h_F \delta
(\|\psi\|_a + \kappa_F^{-1}h_F^{-1}\|\psi\|_0)\lesssim
\varepsilon^{3/2} \kappa_F^{1/2} \| n_F^\top [D_h^2 u_h]_F n_F\|_{0,F}
.
$$
This and the efficiency of $\mu_1^2(T_\pm)$ imply the efficiency of $\mu_2^2(F)$.
The proofs of efficiency of $\mu_3^2(F)$ follows from a more standard
argument \cite{VRposteriori,Verfurth1998}. Indeed, using a standard
$H^2$ conforming face bubble function and integration by parts,
it can be shown that
$$
\mu_3^2(F)
\lesssim
\mu_1^2(T_-)+\mu_1^2(T_+)+\mu_2^2(F)
+\LLL u-u_h\RRR_{\omega_F}^2.
$$
The details are omitted for brevity.
\end{proof}

\begin{remark}
 Parts of the analysis are still valid under the weaker assumption
 that in case $\alpha=1$ the jump $[u_h]_F$ be orthogonal to all
 polynomials from $P_{\ell-1}(F)$ 
 instead of the continuity requirement
 in Assumption~\ref{assumption:main}.
 Our efficiency proof for the error estimator contribution $\xi^2$ 
 will, however, not be robust with respect to $\varepsilon$,
 which is the reason why we assume continuity of $u_h$ in the case
 $\alpha=1$.
\end{remark}

\begin{remark}
 An alternative to our treatment of the tangential-normal part of the 
 Hessian jump is to perform an additional integration by parts
 with the surface gradient on the face $F$ in the formula of
 Lemma~\ref{l:int-by-parts}. This would, however, require stronger
 continuity requirements on the nonconforming finite element space
 in three dimensions than assumed in this paper.
\end{remark}

\section{Some Examples}\label{s:examples}
In this section, we list some nonconforming finite elements to which our theory
applies.
A summary is displayed in Table~\ref{table:elementlist}.
We begin with discussing some classical nonconforming methods.

\begin{example}[Morley element]
The Morley element is based on piecewise quadratic polynomials,
$\ell=2$.
The local degrees of freedom are the evaluations of function averages 
over the $(d-2)$-dimensional hyperfaces and the averages of the normal
derivative over the $(d-1)$-dimensional hyperfaces.
In two dimensions, there exists a higher-order generalization
to enriched $P_\ell(T)$ shape function spaces
by Blum and Rannacher \cite{HR} as an equivalent displacement
method for the classical Hellan--Herrmann--Johnson scheme.
For $\alpha=0$, theses elements satisfy Assumption~\ref{assumption:main},
but not for $\alpha=1$.
A modification with a discrete bilinear form was proposed
by \cite{ModifyMorley,3DmodifiedMorley}.
An a posteriori error analysis of that method
would require different arguments with discrete norms
and is beyond the scope of this article.
\end{example}
      
\begin{example}[Fraeijs de Veubeke (FV) elements]
There are two elements referred to as Fraeijs de Veubeke (FV)
element. Both are two-dimensional first-order elements.
The shape function space of the element FV1 is the sum of 
$P_2(T)$ and three cubic functions
(see \cite{Non1975} for details).
The degrees of freedom are the point evaluations in the vertices
and the edge midpoints and the averages of the normal derivative
along the edges.
If $\alpha=0$, the element FV1 satisfies 
Assumption~\ref{assumption:main} with $\ell=2$.
The element FV2 is based on piecewise cubic polynomials.
The degrees of freedom are the evaluation of the function in the vertices
and in the barycentre and the evaluation of the normal derivative
in two Gauss points of each edge.
This element does not satisfy Assumption~\ref{assumption:main}.
\end{example}

\begin{example}[Specht elements]
The shape function space of the classical triangular Specht element 
is the sum of $P_2(T)$ and three fourth-order polynomials.
The nine local degrees of freedom are the evaluations of the 
function and the gradient in the vertices.
The Specht finite element space is $C^0$ conforming.
As shown in \cite{Specht1ord} the Specht element satisfies 
Assumption~\ref{assumption:main} with $\ell=2$
for $\alpha=0,1$.
A generalization to higher dimensions $d\geq 3$ was proposed in
\cite{Wang-Shi-Xu}, therein named the NZT element.
The second-order version of the Specht element 
proposed in \cite{Shi-Chen-Huang,specht_quadratic,specht2} is based on
a shape function space consisting of $P_3(T)$ plus three quintic
bubble functions. The three additional degrees of freedom compared with
the classical version are the averages of the normal derivatives along
the edges. 
For $\alpha=0,1$,
this element satisfies Assumption~\ref{assumption:main} with $\ell=3$.
\end{example}

\begin{example}[Nilssen--Tai--Winther (NTW) element]
The local shape function space of this triangular element is
the sum of $P_2(T)$ and the three functions resulting from
multiplication of the cubic volume bubble with a barycentric
coordinate. The degrees of freedom are the evaluations in the 
vertices and the edge midpoints and averages of the normal
derivative along the edges.
This $C^0$ conforming element satisfies 
Assumption~\ref{assumption:main} for $\ell=2$ and 
$\alpha=0,1$.
\end{example}

\begin{example}[further simplicial elements]
For $\alpha=0,1$, our error estimator applies to
the Guzm\'an--Leykekhman--Neilan family \cite{Guzm2012} ($d=2,3$, $\ell\geq2$),
the modified Specht element \cite{MT} ($d=2,3$, $\ell\geq 2d-1$),
the Wang-Zu-Zhang element \cite{Wang-Zu-Zhang} ($d=2,3$, $\ell=3$), and
the Chen--Chen--Qiao element  \cite{Chen-Chen-Qiao} ($d=3$, $\ell=3$).
For $\alpha=0$,
the theory applies to the
Hu--Zhang element \cite{Hu&Zhang2019} ($d=2,3$, $\ell=3$) and the
Hu--Tian--Zhang element
\cite{Hu&Zhang2019,HZT} ($d=3$, $\ell\geq3$).
\end{example}

\begin{example}[rectangular elements]
If $\alpha=0$, Assumption~\ref{assumption:main} can be verified with
$\ell=2$ for the 
rectangular Morley element and the incomplete biquadratic
element \cite{Shi86,WangShiXu07} where the degrees of freedom are the evaluations at the 
vertices and the averages of the normal derivative over the edges.
As a negative example, we mention the Adini element \cite{Ciarlet2002,Non1975},
which our theory does not apply to.
\end{example}

\begin{table}
\begin{tabular}{l|l|l|l|l|l}
Element & $d$ & $\ell$ & Ref.\ & $\alpha=0$ & $\alpha=1$ \\
\hline
Morley  & 2,3 & 2      &\cite{Ciarlet2002,Ming2013Minimal}   &yes &no  \\
Blum--Rannacher&2&$\geq3$&\cite{HR}                          &yes &no\\
FV1     & 2   & 2      &\cite{Ciarlet2002,Non1975}           &yes  &no  \\
FV2     & 2   & 3      &\cite{Ciarlet2002,Non1975}           &no  &no  \\
Specht  & 2  & 2,3    &\cite{Specht1ord,specht_quadratic}   &yes &yes \\
NZT     & 2, 3 & 2    &\cite{Specht1ord,Wang-Shi-Xu}   &yes &yes \\
NTW     & 2   & 2      &\cite{NTW}                           &yes &yes \\
Guzm\'an--Leykekhman--Neilan&2,3&$\geq2$&\cite{Guzm2012}&yes&yes\\
modified Specht&2, 3&$\geq 2d-1$&\cite{MT}&yes&yes\\
Wang--Zu--Zhang&2, 3&3&\cite{Wang-Zu-Zhang}&yes&yes\\
Chen--Chen--Qiao&3&3&\cite{Chen-Chen-Qiao}&yes&yes\\
Hu--Zhang&2&$3$&\cite{Hu&Zhang2019} &yes &no \\
Hu--Tian--Zhang&3&$\geq 2$&\cite{Hu&Zhang2019,HZT} &yes &no \\
rectangular Morley    &2& 2   &\cite{WangShiXu07}                 &yes &no \\
incomplete biquadratic &2&2   &\cite{Shi86,WangShiXu07}           &yes &no \\
Adini   & 2 &  3      &\cite{Ciarlet2002}                   &no &no  \\
\end{tabular}
\caption{Applicability of the proposed error estimator to some
existing finite elements.
\label{table:elementlist}}
\end{table}

\section{Numerical experiments}\label{s:num}
In this section we present numerical experiments in two space dimensions.

\subsection{Setup}
The two methods we consider are the Nilssen--Tai--Winther element 
(first-order method, abbreviated by NTW in the diagrams) and 
the modified Specht element from \cite{MT} with polynomial order $\ell=3$
(second-order method, abbreviated by MS in the diagrams).
We use uniform mesh refinement as well as
adaptive meshes generated by the local contributions of our 
error estimator $\eta$ and D\"orfler marking \cite{Doerfler1996}
with marking parameter $1/2$.
Our convergence history plots display the error quantities
in dependence of the number of degrees of freedom (ndof).

\subsection{Example 1: The biharmonic equation}
\begin{figure}
    \begin{center}
	\includegraphics[width=.7\textwidth]{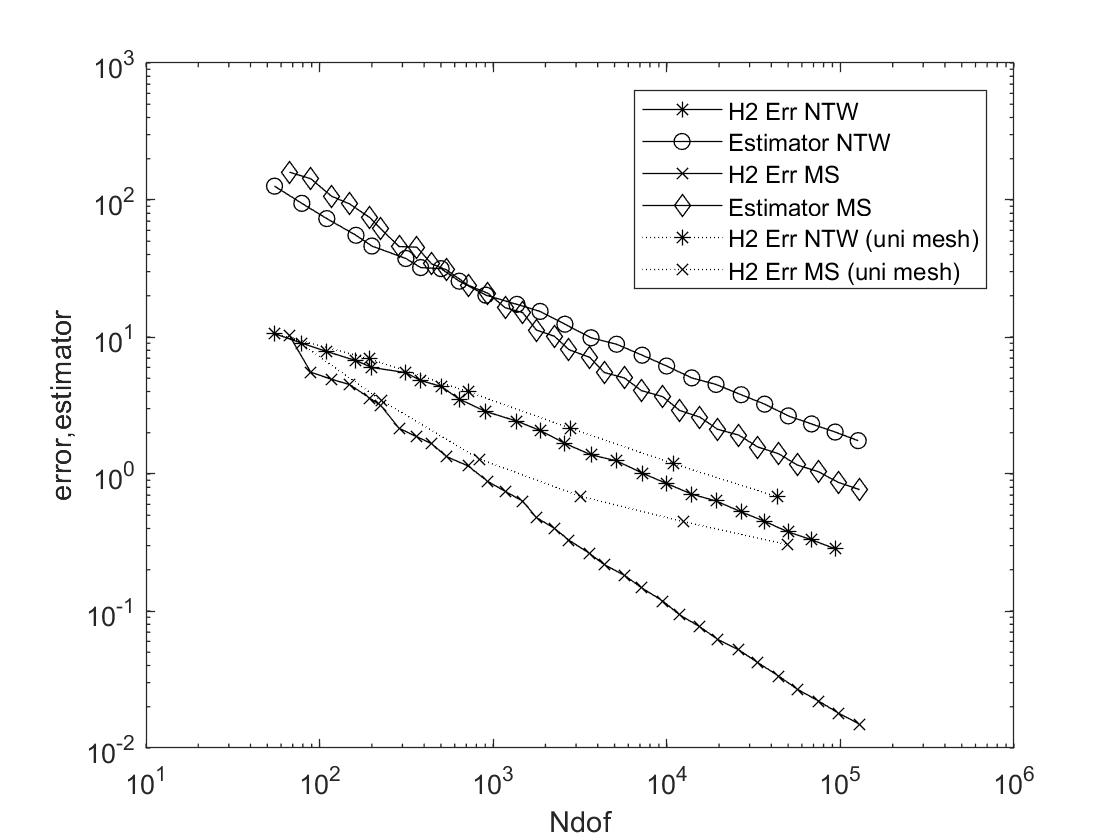}
    \end{center}
	\caption{Convergence history for Example 1.}
	\label{valuebi}
\end{figure}

In the first example we consider the biharmonic equation
with $\varepsilon=1$ and $\alpha=0$.
Let $\Omega= (-1,1)^2\backslash ([0,1]\times [-1,0])$
be the L-shaped domain. The right-hand side $f$ is computed according
to the exact solution \cite{Grisvard1992,BrennerGudiSung2010},
given in polar coordinates by
$$
u(r,\theta)=(r^2\cos^2\theta-1)^2(r^2\sin^2\theta-1)^2r^{1+\alpha}g(\theta),
$$
for $\omega = 3\pi/2$, $\alpha = 0.5444837$,
and the cutoff function
\begin{align*} 
	g(\theta) =\left(\frac{s_-(\omega)}{\alpha-1}-\frac{s_+(\omega)}{\alpha+1} \right)(c_-(\theta)-c_+(\theta))
	-\left(\frac{s_-(\theta)}{\alpha-1}-\frac{s_+(\theta)}{\alpha+1} \right)(c_-(\omega)-c_+(\omega)).
\end{align*}
with the abbreviations $s_\pm(z)=\sin((\alpha\pm1)z)$
and $c_\pm(z)=\cos((\alpha\pm1)z)$.
Figure~\ref{valuebi} shows the convergence history of 
the absolute $\|\cdot\|_{a,h}$ norm errors and the error estimators $\eta$.
We observe the proven equivalence of the error and the error estimator.
Furthermore, we observe that the adaptive method converges with the optimal rate,
while on uniform meshes the methods converge with the expected suboptimal rate.

\subsection{Example 2: Singularly perturbed problem on the square}

\begin{figure}
    \begin{center}
	\includegraphics[width=.7\textwidth]{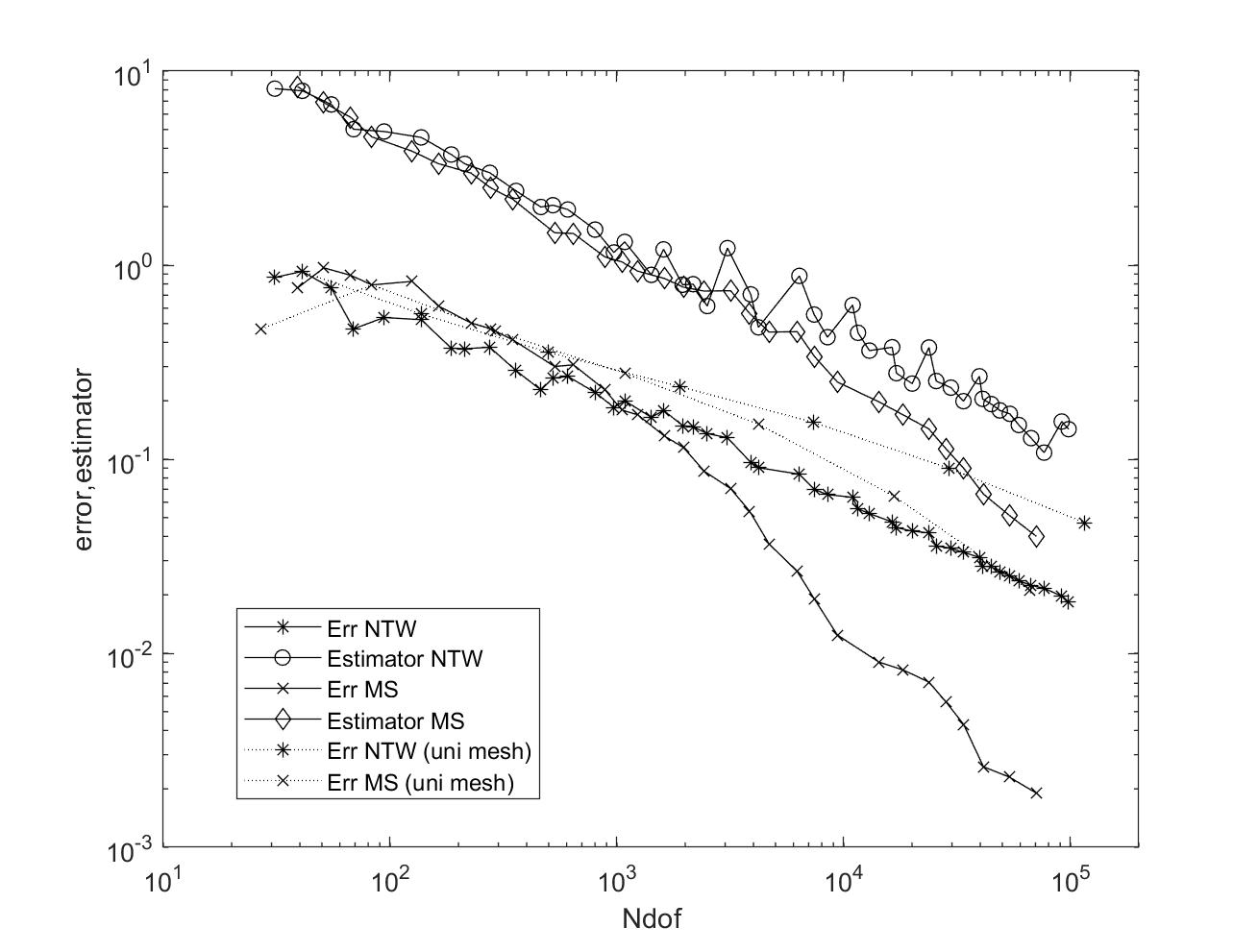}
	\end{center}
	\caption{Convergence history for Example 2.}
	\label{valuebl}
\end{figure}

In this example, we consider the unit square $\Omega=(0,1)^2$
and parameters $\alpha=1$, $\varepsilon=10^{-2}$.
The exact solution
$$
u = (\sin \pi x-\pi \varepsilon\frac{\cosh\frac{1}{2\varepsilon} -\cosh\frac{2x-1}{2\varepsilon}}{\sinh\frac{1}{2\varepsilon}} )(\sin \pi y-\pi \varepsilon\frac{\cosh\frac{1}{2\varepsilon} -\cosh\frac{2y-1}{2\varepsilon}}{\sinh\frac{1}{2\varepsilon}} )
$$
has a boundary layer due to the incompatibility of the boundary
condition with the limiting solution
$$
\lim_{\varepsilon\rightarrow 0}u = \sin\pi x\sin \pi y.
$$
Figure~\ref{valuebl} displays the convergence history
the absolute $\|\cdot\|_{a,h}$ norm errors and the error estimators $\eta$.
Due to the boundary layer, the asymptotic convergence regime
starts after approximately $2\,000$ degrees of freedom only.
After that, the adaptive methods show the optimal convergence rate, which 
indicates that the layer has been resolved by the adaptive method.
The optimal convergence rate with uniform refinement is observed 
starting from about $10^4$ degrees of freedom.

\subsection{Example 3: Singularly perturbed problem on the L-shaped domain}
\begin{figure}
    \begin{center}
	\includegraphics[width=.7\textwidth]{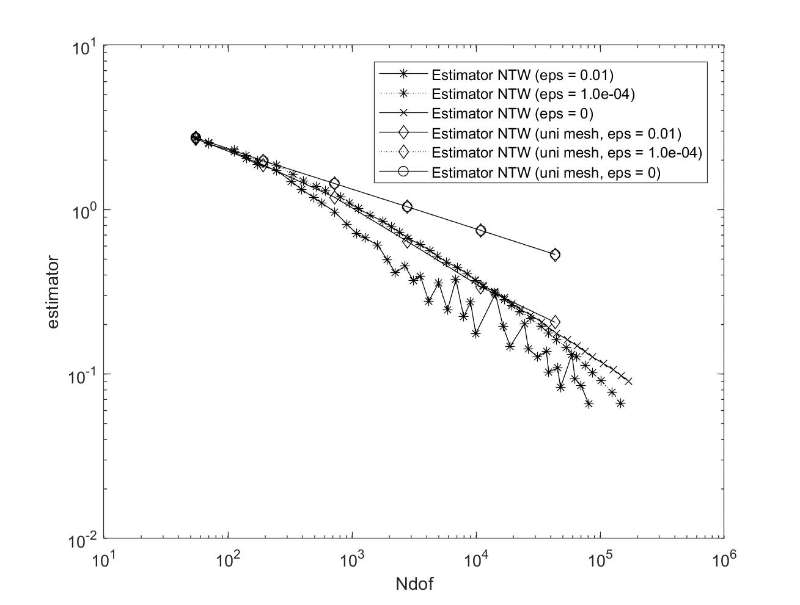}
	\end{center}
	\caption{Values of the error and estimator for Example 3 with the NTW element. }
	\label{valuebl2}
\end{figure}

\begin{figure}
    \begin{center}
	\includegraphics[width=.7\textwidth]{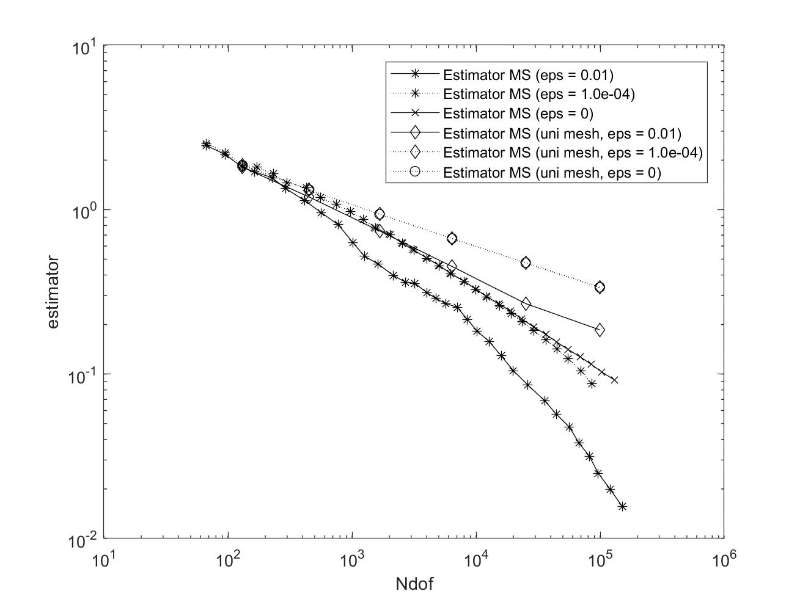}
	\end{center}
	\caption{Values of the error and estimator for Example 3 with the MS element. }
	\label{valuebl3}
\end{figure}

\begin{figure}
    \begin{center}
	\includegraphics[width=.4\textwidth]{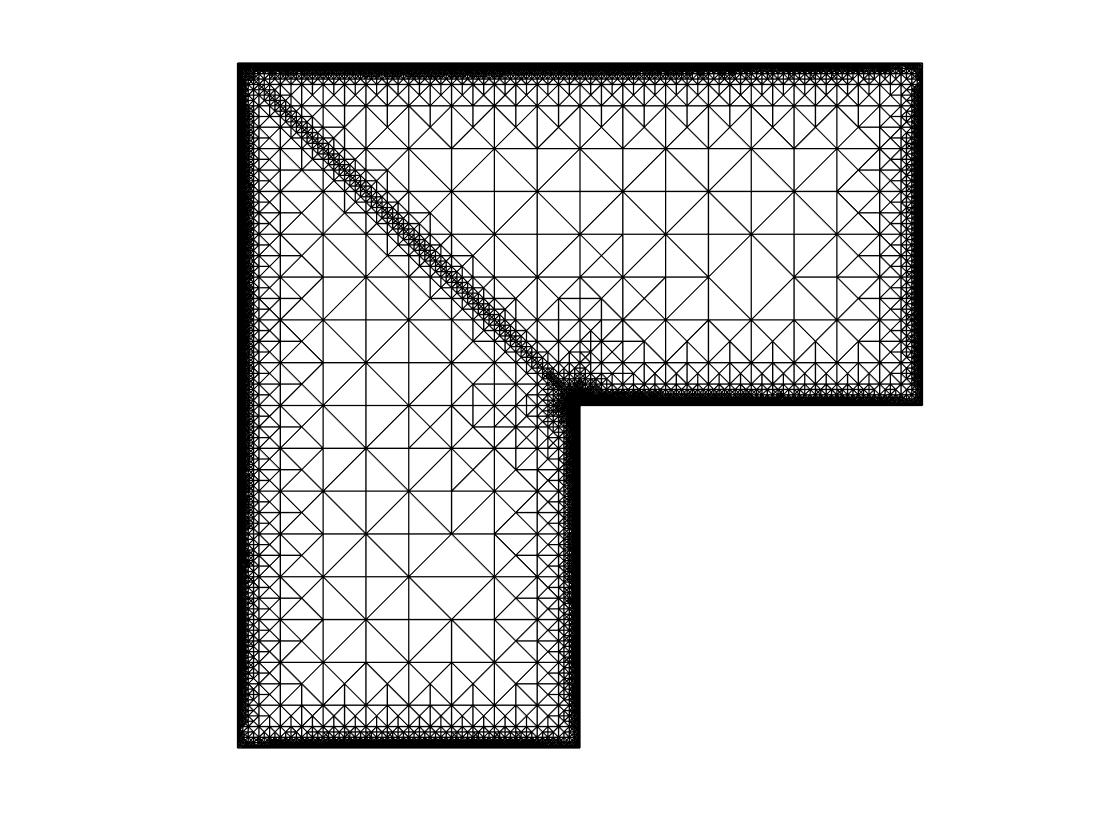}
	\includegraphics[width=.4\textwidth]{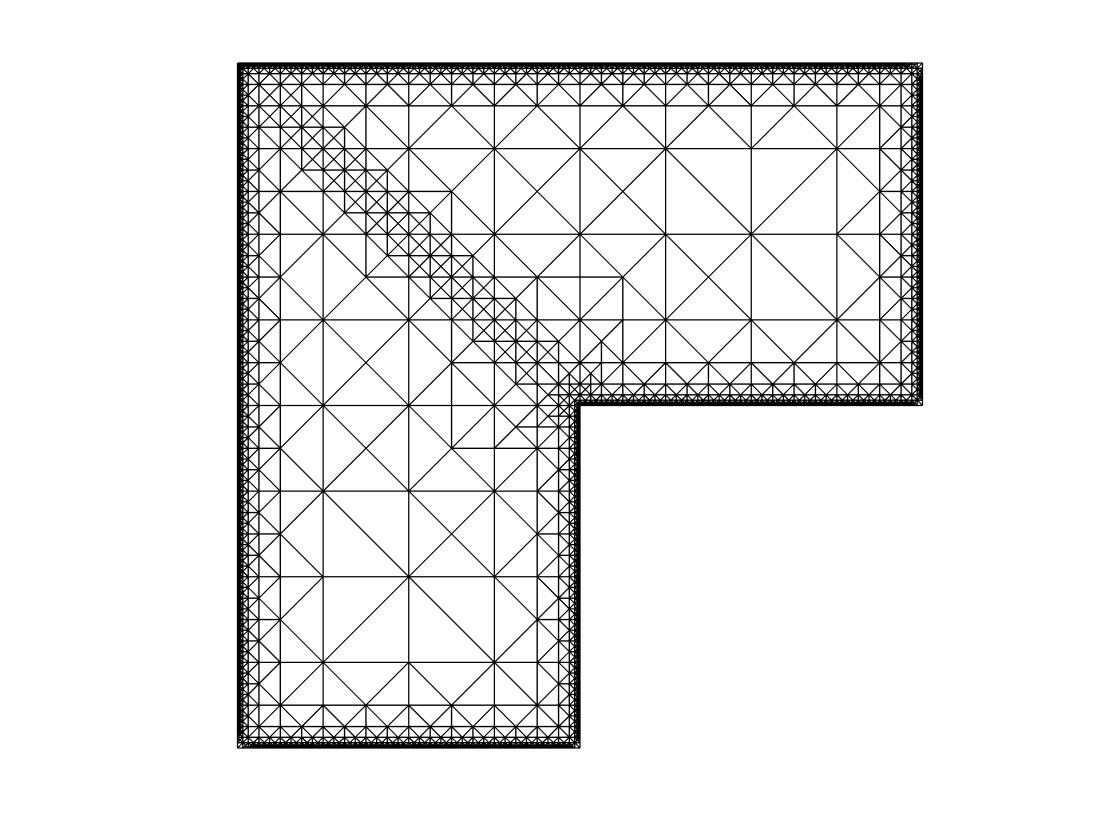}
	\end{center}
	\caption{Adaptive mesh with 120\,577 (left) and 130,617 (right) degrees of freedom generated by MS element
	       when $\varepsilon=1/100$ and $\varepsilon=0$, respectively, in Example 3.}
	\label{mesh}
\end{figure}

We consider the L-shaped domain, $\alpha=1$ and 
$\varepsilon \in \{1/100, 10^{-4},0\}$
and the the right hand side $f(x,y)= (|x+y|)^{-1/3}$.
In this example, no exact solution is known.
Figure~\ref{valuebl2} and \ref{valuebl3} display the convergence history
the error estimators $\eta$.
For the choice $\varepsilon=1/100$, both methods are observed
to converge at the optimal rate.
Two meshes with $120\,577$ and $130\,617$ degrees of freedom generated with the MS element
are displayed in Figure~\ref{mesh}.
For the more challenging choice $\varepsilon=10^{-4}$, all methods
show only first-order rates within the range of numbers of degrees
of freedom under consideration, which indicates that the boundary
layer has not been resolved, yet.
For the limit case $\varepsilon = 0$, we first note that due to 
the potentially incompatible boundary conditions, the exact 
solution will usually not satisfy the boundary condition for the 
normal derivative. Accordingly,
in the numerical experiment, we observe convergence of the error estimator
at the rate of a first-order method.
This is due to the fact that the normal-derivative degrees of freedom
are set to $0$ on the boundary and thereby basis functions are sorted 
out that would grant a higher approximation order. On the other hand,
we observe robustness of the method as well as reliability and 
efficiency of the error estimator in the limit case
(which is not difficult to prove with the techniques established above).

\appendix
\section{Generalized Hsieh--Clough--Tocher elements}\label{ax:HCT}
For $d=2$, the interpolation operator is defined on the HCT splits \cite{HCT}
and for $d=3$ on the Worsey--Farin splits \cite{WS}.
The degrees of freedom for interpolation  in the lowest order case
(piecewise cubic) on both splits are the function value and the gradient value 
on the vertices and the directional derivatives at the midpoints
in $d-1$ directions normal to the edge.
The high order case can be found in \cite{HighorderHCT,WSNeilan}.
For convenient reading, we briefly describe the construction here.

Given a simplex $K\subseteq \mathbb R^d$ triangulated by a subdivision
$\mathcal T_K$, we denote by $\mathcal{P}_{\ell}(\mathcal T_K)$ the space
of $C^1(K)$ functions that are piecewise polynomials of degree $\leq \ell$
with respect to $\mathcal T_K$.
We furthermore denote
$$
 \tilde{\mathcal{P}}_{\ell}(\mathcal T_K)
 :=\{v\in \mathcal P_\ell(\mathcal T_K)
                 |\nabla v|_{\partial K}=0\}
\quad\text{and}\quad
\tilde{\mathcal{P}}_{\ell,0}(\mathcal T_K)
:=\{v\in\tilde{\mathcal{P}}_{\ell}(\mathcal T_K)| v_{\partial K}=0\}.
$$

\subsection{High order HCT element in 2D}
Figure \ref{f:wfhct} shows an HCT split. The three vertices of the triangle $K$
are denoted as $a_1, a_2, a_3$ and the edge opposite to $a_i$ is
denoted as $e_i$, $i = 1, 2, 3$.
Let as $c_K$ denote a given interior point of $K$.
The sub-triangulation $\mathcal T_{K,c_K}^{\mathrm{HCT}}$
is formed by dividing $K$ in three sub-triangles $K_i, i = 1, 2, 3$, 
where $K_i$ is the convex hull of $c$ and $e_i$.
The three interior edges $g_i$, $i=1,2,3$, are the 
line segments connecting $a_i$ and $c$.
The high order ($\ell\geq 4$) HCT element is determined by the following
degrees of freedom.
A function $v\in \mathcal P_\ell(\mathcal T_{K,c_K}^{\mathrm{HCT}})$
is uniquely determined by
\begin{align*}
     &v(a_i),\nabla v(a_i),\quad  v(c),\nabla v(c), &&\\
         &\fint_{e_i}v q_{\ell-4}\mathrm{d}s & &\text{ for all } q_{\ell-4}\in
P_{\ell-4}(e_i),\\
         &\fint_{e_i}\partial_{n_{e_i}}v q_{\ell-3}\mathrm{d}s, & &\text{ for all }q_{\ell-3}\in
P_{\ell-3}(e_i),\\
         &\fint_{g_i}\partial_{n_{g_i}}v q_{\ell-5}\mathrm{d}s,~\fint_{g_i}v
q_{\ell-5}\mathrm{d}s  & &\text{ for all } q_{\ell-5}\in P_{\ell-5}(g_i),\\
         & \fint_{ K_i}vq_{\ell-6}\mathrm{d}s & &\text{ for all } q_{\ell-6}\in
P_{\ell-6}(K_i),
\end{align*}
for $i=1,2,3$.
For proofs we refer to \cite{DouglasDupontPercellScott1979}.

\subsection{High order element in 3D on Worsey-Farin splits}
\begin{figure}[tb] \setlength\unitlength{4pt}		
\begin{tikzpicture}[scale=2.0]
	\draw [color=white] (1,-0.5)--(1,0);
	\draw (0,0)--(2,0);
	\draw (1,1.5)--(0,0);
	\draw (1,1.5)--(2,0);
	\draw (1,0.5)--(0,0);
	\draw (1,0.5)--(2,0);
	\draw (1,0.5)--(1,1.5);
	\node [below] at (1,0.5) {$c$};
	\fill (1,0.5) circle (1.0pt); 
	\fill (0,0) circle (1.0pt); 
	\node [left] at (0,0) {$a_1$};
	\node [right] at (2,0) {$a_2$};
	\fill (2,0) circle (1.0pt); 
	\node [above] at (1,1.5) {$a_3$};
	\fill (1,1.5) circle (1.0pt); 
	\node [left] at (0.9,0.6) {$K_2$};
	\node [right] at (1.1,0.6) {$K_1$};
	\node [below] at (1,0.3) {$K_3$};
	\node [left] at (0.5,0.75) {$e_2$};
	\node [right] at (1.5,0.75) {$e_1$};
	\node [below] at (1,0) {$e_3$};
	\node [left] at (1,0.9) {$g_3$};
	\node [above] at (1.5,0.25) {$g_2$};
	\node [below] at (0.6,0.3) {$g_1$};
\end{tikzpicture}
\qquad
	\begin{tikzpicture}[scale=1.8]
		\draw (0,0)--(1,1.8);
		\draw  (1,1.8)--(3,0);
		\draw [dashed]  (0,0)--(3,0);
		\draw  (0,0) -- (1.5,-0.5);  
		\draw  (1.5,-0.5)--(3,0);
		\draw  (1.5,-0.5)--(1,1.8);
		\fill (0,0) circle (1.5pt); 
		\fill (1,1.8) circle (1.5pt); 
		\fill (3,0) circle (1.5pt); 
		\fill (1.5,-0.5) circle (1.5pt);
		\fill (0.8,0.6) circle(1.5pt); 
		\node [left] at (0.8,0.6) {$c_1$};
		\fill (1.5,0.5) circle(1.5pt);
		\node [below] at (1.5,0.5) {$c_T$};
		\draw [dashed] (0.8,0.6)--(1.5,0.5);
		\draw [dashed] (1.5,-0.5)--(1.5,0.5);
		\draw [dashed] (1,1.8)--(1.5,0.5);
		\draw [dashed] (0,0)--(1.5,0.5);
		\draw (0.8,0.6) -- (0,0);
		\draw (0.8,0.6) -- (1.5,-0.5);
		\draw (0.8,0.6) -- (1,1.8);
		\draw  [dashed] (1.5,0.5) -- (3,0);
		\node [above] at (1,1.8) {$a_2$};
		\node [left] at (0,0) {$a_3$};
	   \node [below] at (1.5,-0.5) {$a_4$};
	   	\node [right=] at (3,0) {$a_1$};
	\end{tikzpicture}

\caption{Hsieh--Clough--Tocher split (left) and  Worsey--Farin split (right).
        \label{f:wfhct}}
\end{figure}
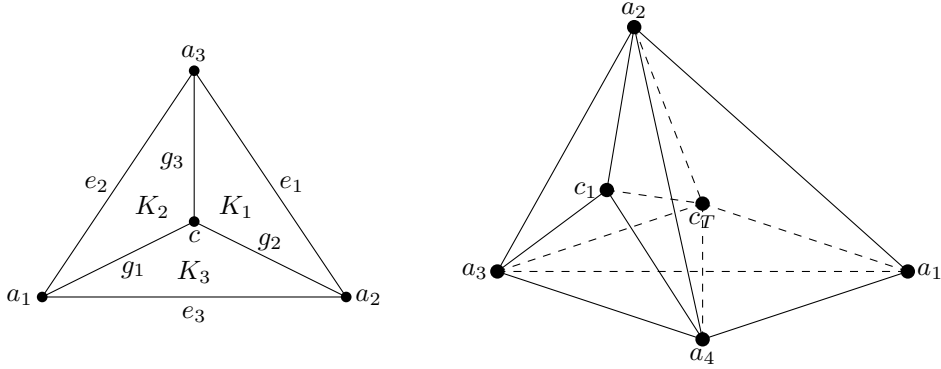

The four vertices of the simplex $K$ are denoted as $a_1,a_2,a_3,a_4.$
Let $c_K$ denote the midpoint of $K$.
Let $F_i$ denote the face opposite to $a_i$.
If $F_i$ is an interior face of the triangulation containing $K$,
it is shared by two element $K$ and $K'$,
and $c_i\in F_i$ is defined as 
$c_i=\overline{c_Kc_{K'}}\cap F_i$,
i.e., the intersection of the line connecting
$c_K$ and $c_{K'}$, and $F_i$.
The Worsey--Farin split,
displayed in Figure~\ref{f:wfhct}, consists of the 12 simplices 
$[c_K,a_i,a_j,c_k],~1\leq i\neq j\neq k\leq 4$ 
and is denoted by $\mathcal T_K^{\mathrm{WF}}$.
The high order ($\ell\geq 4$) Worsey--Farin element is determined by the following
degrees of freedom.
A function $v\in \mathcal P_\ell(\mathcal  T_K^{\mathrm{WF}})$
is uniquely determined by
\begin{align*}
	&v(a_m), \nabla v(a_m),&&\\
	&\fint_{e_i}vq_{\ell-4}\mathrm{d}s 
	  &&\text{ for all } q_{\ell-4}\in P_{\ell-4}(e_i),\\
	&\fint_{e_i}\partial_{n_{e,j}}vq_{\ell-3}\mathrm{d}s 
	  &&\text{ for all } q_{\ell-3}\in P_{\ell-3}(e_i),\\
    &\fint_{F_m}(\nabla_{F_m}v\cdot g)dS 
      &&\text{ for all } g\in \tilde{\mathcal{P}}_{\ell,0}(\mathcal T_{F_m,c_m}^{\mathrm{HCT}}),\\
    &\fint_{F_m}(\partial_{n_{F_m}} v)g\mathrm{d}S 
      &&\text{ for all } g\in \tilde{\mathcal{P}}_{\ell-1}(\mathcal T_{F_m,c_m}^{\mathrm{HCT}}),\\
    &\fint_T \nabla v\cdot g\mathrm{d}x 
      &&\text{ for all } g\in \nabla \tilde{\mathcal{P}}_{\ell,0}(\mathcal T_K^{\mathrm WF}).
\end{align*}
for all $i=1,\dots,6$, $m=1,\dots,4$, and $j=1,2$.
Here $\nabla_F$ is the surface gradient on $F$.
The six edges of $K$ are denoted by $e_1,\dots,e_6$,
and $n_{e_i,j},~j=1,2$ are two linear 
independent vectors orthogonal to the edge $e_i$. 
For proofs we refer to \cite{WSNeilan}.

\section{Proof of Lemma~\ref{l:localization}}\label{ax:localizationproof}

Without loss of generality we consider the most relevant case $\alpha=1$.
We consider a $H^2$-conforming partition of unity with respect to
the triangulation $\mathcal T$ and its vertices $\mathcal N$.
This family $(\varphi_y)_{y\in\mathcal N}$
consists of continuous and bounded functions and has the property
$\sum_{y\in\mathcal N} \varphi_y= 1$ in $\Omega$
and $\operatorname{supp}(\varphi_j) = \overline \omega_y$
where $\omega_y$ is the vertex patch of diameter $h_y$.
The partial derivatives satisfy
\begin{equation}\label{e:scaling_varphi}
 \max_{\omega_y}| \partial^\beta \varphi_y|
 \lesssim
 h_y^{-|\beta|}
\end{equation}
for any multi-index $\beta$ of length $|\beta|\leq 2$.
Such partition is provided by the basis functions of a $H^2$-conforming
finite element space with respect to $\mathcal T$ corresponding to the 
point evaluation in the vertices.
Since any function of the format
$$
  v = \sum_{y\in\mathcal N}  v_y \varphi_y
$$
with functions $v_y\in V(\omega_y)$ belongs to $V$,
and since the functions $\varphi_j$ form a partition of unity,
we can write
$$
\min_{v\in V} \LLL u_h - v \RRR_{}^2
=
\min_{(v_y)_y\in V(\omega_y)_y}
\LLL \sum_{y\in\mathcal N} \varphi_y(u_h - v_y) \RRR_{}^2
\lesssim
\sum_{y\in\mathcal N}
\min_{v_y\in V(\omega_y)}
\LLL  \varphi_y(u_h - v_y) \RRR_{}^2
$$
due to the bounded overlap of vertex patches.
We abbreviate $w:=u_h - v_y$.
From the product rule and the scaling \eqref{e:scaling_varphi}
we infer
$$
\LLL  \varphi_y w \RRR_{}^2
\lesssim
(h_y^{-2}+\varepsilon^2h_y^{-4}) \|w\|_{0,\omega_y}^2 
+
(1+\varepsilon^2h_y^{-2}) \|\nabla_h w \|_{0,\omega_y}^2
+
\varepsilon^2 \|D^2_h w\|_{0,\omega_y}^2.
$$
From $(1+\varepsilon^2/h_y^2)\leq 2 \kappa_y^{-2}$
we deduce
$$
\LLL  \varphi_y w \RRR_{}^2
\lesssim
h_y^{-2}\kappa_y^{-2} \|w\|_{0,\omega_y}^2 
+
\kappa_y^{-2} \|\nabla_h w\|_{0,\omega_y}^2
+
\varepsilon^2 \|D^2_h w\|_{0,\omega_y}^2.
$$
The combination of the above arguments concludes the proof.

\bibliographystyle{plain}
\bibliography{posterior}

\end{document}